\documentclass[12pt]{amsart}
\usepackage{amssymb}

\begin{document}
\numberwithin{equation}{section}

\def\Label#1{\label{#1}}

\def\1#1{\ov{#1}}
\def\2#1{\widetilde{#1}}
\def\3#1{\mathcal{#1}}
\def\4#1{\widehat{#1}}

\def\s{s}
\def\k{\kappa}
\def\ov{\overline}
\def\span{\text{\rm span}}
\def\tr{\text{\rm tr}}
\def\GL{{\sf GL}}
\def\xo {{x_0}}
\def\Rk{\text{\rm Rk\,}}
\def\sg{\sigma}
\def \emxy{E_{(M,M')}(X,Y)}
\def \semxy{\scrE_{(M,M')}(X,Y)}
\def \jkxy {J^k(X,Y)}
\def \gkxy {G^k(X,Y)}
\def \exy {E(X,Y)}
\def \sexy{\scrE(X,Y)}
\def \hn {holomorphically nondegenerate}
\def\hyp{hypersurface}
\def\prt#1{{\partial \over\partial #1}}
\def\det{{\text{\rm det}}}
\def\wob{{w\over B(z)}}
\def\co{\chi_1}
\def\po{p_0}
\def\fb {\bar f}
\def\gb {\bar g}
\def\Fb {\ov F}
\def\Gb {\ov G}
\def\Hb {\ov H}
\def\zb {\bar z}
\def\wb {\bar w}
\def \qb {\bar Q}
\def \t {\tau}
\def\z{\chi}
\def\w{\tau}
\def\Z{\zeta}

\def \T {\theta}
\def \Th {\Theta}
\def \L {\Lambda}
\def\b{\beta}
\def\a{\alpha}
\def\o{\omega}
\def\l{\lambda}

\def \im{\text{\rm Im }}
\def \re{\text{\rm Re }}
\def \Char{\text{\rm Char }}
\def \supp{\text{\rm supp }}
\def \codim{\text{\rm codim }}
\def \Ht{\text{\rm ht }}
\def \Dt{\text{\rm dt }}
\def \hO{\widehat{\mathcal O}}
\def \cl{\text{\rm cl }}
\def \bR{\mathbb R}
\def \bC{\mathbb C}
\def \bP{\mathbb P}
\def \C{\mathbb C}
\def \bL{\mathbb L}
\def \bZ{\mathbb Z}
\def \bN{\mathbb N}
\def \scrF{\mathcal F}
\def \scrK{\mathcal K}
\def \scrM{\mathcal M}
\def \cR{\mathcal R}
\def \scrJ{\mathcal J}
\def \scrA{\mathcal A}
\def \scrO{\mathcal O}
\def \scrV{\mathcal V}
\def \scrL{\mathcal L}
\def \scrE{\mathcal E}
\def \hol{\text{\rm hol}}
\def \aut{\text{\rm aut}}
\def \Aut{\text{\rm Aut}}
\def \J{\text{\rm Jac}}
\def\jet#1#2{J^{#1}_{#2}}
\def\gp#1{G^{#1}}
\def\gpo{\gp {2k_0}_0}
\def\emmp {\scrF(M,p;M',p')}
\def\rk{\text{\rm rk}}
\def\Orb{\text{\rm Orb\,}}
\def\Exp{\text{\rm Exp\,}}
\def\ess{\text{\rm Ess\,}}
\def\mult{\text{\rm mult\,}}
\def\Jac{\text{\rm Jac\,}}
\def\Span{\text{\rm span\,}}
\def\d{\partial}
\def\D{\3J}
\def\pr{{\rm pr}}
\def\dbl{[\![}
\def\dbr{]\!]}
\def\nl{|\!|}
\def\nr{|\!|}

\def \depth{\text{\rm depth\,}}
\def \D{\text{\rm Der}\,}
\def \Rk{\text{\rm Rk}\,}
\def \ima{\text{\rm im}\,}
\def \vfi{\varphi}

\title[Analyticity of smooth CR mappings of generic
submanifolds] {Analyticity of smooth CR
mappings of generic
submanifolds}
\author[P. Ebenfelt and L. P. Rothschild]{Peter Ebenfelt and
Linda P. Rothschild} \footnotetext{{\rm The first author is
supported in part by DMS-0401215. The second author is
supported in part by DMS-0400880.\newline}}
\address{ Department of Mathematics, University of California
at San Diego, La Jolla, CA 92093-0112, USA}
\email{pebenfel@math.ucsd.edu, lrothschild@ucsd.edu }



\abstract We consider a smooth CR mapping $f$ from a real-analytic generic submanifold $M$ in $\bC^N$ into $\bC^N$.  For $M$  of finite type and essentially finite at a point $p\in M$, and 
$f$ formally finite at $p$, we give a necessary and sufficient condition
for $f$ to extend as a holomorphic mapping in some neighborhood of $p$.  In a similar vein, we consider a formal holomorphic mapping $H$ and give a necessary and sufficient condition for $H$ to be convergent.\endabstract

\dedicatory{
Dedicated to Salah Baouendi for his seventieth birthday.
}

\newtheorem{Thm}{Theorem}[section]
\newtheorem{Def}[Thm]{Definition}
\newtheorem{Cor}[Thm]{Corollary}
\newtheorem{Pro}[Thm]{Proposition}
\newtheorem{Lem}[Thm]{Lemma}
\newtheorem{Rem}[Thm]{Remark}
\newtheorem{Ex}[Thm]{Example}

\maketitle
\section{Introduction}

In this paper, we consider a smooth CR mapping $f$ from a real-analytic generic submanifold $M$ in $\bC^N$ into $\bC^N$.  Assuming that $M$ is of finite type and essentially finite at a point $p\in M$, and that
$f$ is formally finite at $p$ (see below), we give a necessary and sufficient condition
for $f$ to extend as a holomorphic mapping (Theorem \ref{t:ext1}) in some neighborhood of $p$ (or
equivalently to be real-analytic near $p$ in $M$).  In a similar vein, we consider a formal holomorphic mapping $H$ and give a necessary and sufficient condition for $H$ to be convergent (Theorem  \ref{t:ext2}).

Before stating the main results, we shall recall some definitions.
Let $M$ be a  real-analytic submanifold of codimension $d$ in
$\bC^N$. Recall that $M$ is said to be {\it generic} if $M$ is
defined locally near any point $p\in M$ by defining equations
$\rho_1(Z,\bar Z)=\ldots=\rho_d(Z,\bar Z)=0$ such that
$\partial\rho_1\wedge\ldots\wedge\partial\rho_d\neq 0$ along $M$.  A
generic submanifold $M$ is said to be of {\it finite type} at $p\in
M$ (in the sense of Kohn \cite{Kohn} and Bloom-Graham \cite{BG}) if
the (complex) Lie algebra $\frak g_M$ generated by all smooth
$(1,0)$ and $(0,1)$ vector fields tangent to $M$ satisfies $\frak
g_M(p_0)=\bC T_{p}M$, where $\C T_{p} M$ is the complexified tangent
space to $M$. For the definition of {\it essentially finite}, the
reader is referred to \cite{BER99a}; see also Section
\ref{s:essvar} for an equivalent formulation in a slightly more
general setting.

A ($C^\infty$) smooth mapping $f\colon
M\to \bC^N$ is called CR if the tangent mapping $df$
sends the CR bundle
$T^{0,1}M$ into $T^{0,1}\bC^N$. In particular, the
restriction to $M$ of a holomorphic mapping $H\colon
U\to \bC^N$, where $U$ is some open neighborhood of $M$, is CR.  To define the notion of {\it formally finite} at a point $p\in M$,  we may assume, without loss of generality,  that $p=f(p)=0$. If  $f\colon
(M,0)\to (\bC^N,0)$ is a germ of a smooth CR mapping, then one may
associate to $f$ a formal mapping $H\colon (\bC^N,0)\to (\bC^N,0)$
as follows. Let $x$ be a local coordinate on $M$ near $0$ and
$x\mapsto Z(x)$ the local embedding of $M$ into $\bC^N$ near $0$.
Then, there is a unique formal mapping $H$ such that the Taylor
series of $f(x)$ at $0$ equals $H(Z(x))$ (see e.g.\ \cite{BER99a},
Proposition 1.7.14).
Recall
that a formal holomorphic (or simply formal) mapping
$H\colon (\bC^N,0)\to (\bC^{N},0)$ is an
$N$-tuple $H=(H_1,\ldots, H_N)$ with $H_j\in
\bC[\![Z_1,\ldots, Z_N]\!]=\bC[\![Z]\!]$ such that each
component
$H_j$ has no constant term. We shall use the notation
$I(H_1,\ldots, H_N)$ (or simply
$I(H)$) for the ideal in $\bC[\![Z]\!]$ generated by the
components $H_1(Z),\ldots, H_N(Z)$. The formal mapping $H$ is
said to be {\it finite} if the ideal $I(H)$ is of finite
codimension in $\bC[\![Z]\!]$, i.e.\ if
$\bC[\![Z]\!]/I(H)$ is a finite dimensional vector space over
$\bC$.  We shall say that the CR mapping $f$ is {\it
formally finite} if the associated formal mapping $H$ is finite. The reader is referred to \cite{BER99a} for further basic notions and properties of generic submanifolds in $\bC^N$ and their mappings.

We may now state our first main result.

\begin{Thm}\Label{t:ext1} Let $M\subset\bC^N$ be a real-analytic
generic submanifold of dimension $m$ that is essentially finite
and of finite type at $0$. Let $f\colon (M,0)\to (\bC^N,0)$ be a
smooth formally finite CR mapping. The following are equivalent.
\medskip

{\rm (i)} There exists an irreducible real-analytic subvariety $\tilde X\subset
\bC^N$ of dimension $m$ at $0$ such that $f(M)\subset \tilde
X$ as germs at $0$.\medskip

{\rm (ii)} $f$ is real-analytic in a neighborhood of $0$.

\end{Thm}

Our second result concerns the convergence of a formal holomorphic mapping sending a real-analytic
generic submanifold into a real-analytic subvariety of the same dimension. Recall that a formal
mapping $H\colon(\bC^N,0)\to (\bC^N,0)$ is said to send a real-analytic generic submanifold $M$ through
$0$ in $\bC^N$ into a real-analytic subvariety $\tilde X$ through $0$ in $\bC^N$, denoted
$H(M)\subset \tilde X$, if $\sigma(H(Z(x)),\overline{H(Z(x))})\equiv 0$ as a power series in
$x$, where $x$ is a local coordinate on $M$ near $0$, $x\mapsto Z(x)$ the local embedding of
$M$ into $\bC^N$ near $0$,
$\sigma(Z,\bar Z)=(\sigma_1(Z,\bar Z),\ldots, \sigma_k(Z,\bar Z))$, and $\sigma_1(Z,\bar Z),
\ldots, \sigma_k(Z,\bar Z)$ generate the ideal of germs at $ 0$ of real-analytic functions vanishing on $\tilde X$.

\begin{Thm}\Label{t:ext2} Let $M\subset\bC^N$ be a real-analytic
generic submanifold of dimension $m$ that is essentially finite
and of finite type at $0$. Let $H\colon (\bC^N,0)\to (\bC^N,0)$ be a finite formal holomorphic mapping. The following are equivalent.
\medskip

{\rm (i)} There exists an irreducible real-analytic subvariety $\tilde X\subset
\bC^N$ of dimension $m$ at $0$ such that $H(M)\subset \tilde
X$.\medskip

{\rm (ii)} $H$ is convergent in a neighborhood of $0$.

\end{Thm}

The proofs of (i)$\implies$(ii) in Theorems \ref{t:ext1} and
\ref{t:ext2}  rest on general criteria for analyticity of smooth CR
mappings and convergence of formal mappings given in \cite{MMZ02}
and \cite{Mirg}, respectively, and a geometric result given in Lemma
\ref{t:essvar} below. A special case of (i) $\implies$ (ii) in Theorem \ref{t:ext1}, in
which the additional hypothesis that the target $\tilde X$ is a
real-analytic generic submanifold $\widetilde M$ of dimension $m$ is
imposed, can be proved by using known results as follows. In
\cite{Me1}, it was shown that $f$ is analytic provided that either
$f$ is CR transversal to $\widetilde M$ at $0$, or $\widetilde M$ is
essentially finite at $0$. The desired conclusion then follows from
a result in \cite{ER05}, where the CR transversality of $f$ was
proved under the hypotheses given. One of the difficulties in the
general case addressed here is the fact that the target $\tilde X$
need not be smooth at $0$ and, consequently, there is no notion of
tranversality of the mapping. This is overcome by showing directly
that the target must satisfy a generalized essential finiteness
condition (see Lemma \ref{t:essvar}), and applying the result from
\cite{MMZ02} mentioned above.

The case where $M$ and $\tilde X$ are real-analytic (non-singular)
hypersurfaces has a long history, beginning with the work of Lewy
\cite{Le} and Pinchuk  \cite {Pi}. There are also many subsequent
results implying analyticity of a smooth CR mapping when the target
$\tilde X$ is a real-analytic generic submanifold $\widetilde M$ of
$\bC^{N'}$ under various hypotheses on $\widetilde M$ and $f$. We
mention here only, in addition to \cite{MMZ02}, the works \cite{DW},
\cite{Han1}, \cite{BJT}, \cite{BRgerms}, \cite{DF}, \cite{Forst}, \cite{Pu},
\cite{BHR}, \cite{Hu96}, \cite{CPS}, \cite{D01}, \cite{MMZ03b} and
refer the reader to the bibliographies of these for further
references. Previous results on convergence of formal mappings were
given e.g.\ in the papers \cite{BER00}, \cite{Mirh}, \cite{BMR02},
\cite{Mirg}, \cite{Mi02b}, \cite{MMZ03a}.

\section{The essential variety of a real-analytic subvariety in $\bC^N$}
\Label{s:essvar}

We begin by defining the notion of essential finiteness for a real-analytic
subvariety $X$ through $0$ in $\bC^N$.  Let  $C^\omega_0$ denote the ring
of germs of real-valued real-analytic functions at $0$ and $I_\bR(X)$ be the ideal in
$C^\omega_0$ of functions vanishing on $X$. Let $\sigma(Z,\bar
Z):=(\sigma_1(Z,\bar Z),\ldots, \sigma_d(Z,\bar Z))$ be
(representatives of) generators of $I_\bR(X)$. We may assume that $\sigma(Z,\zeta)$ is
defined in $B\times B$, where $B$ is a sufficiently small open ball in
$\bC^N$ centered at
the origin. Define the Segre variety
$\Sigma_p\subset B$ of
$X$ at
$p$, for
$p\in B  $, by the holomorphic equations $\sigma(Z,\bar p)=0$. Note that
$\Sigma_0$ is a complex analytic variety through $0$ and, by the
reality of the functions $\sigma(Z,\bar Z)$, the variety $\Sigma_p$ passes through $0$ for every $p\in \Sigma_0$.
Let $U\subset B$ be an open neighborhood of $0$ and define
\begin{equation}\Label{e:essvar}
E^U_0:=\bigcap_{p\in \Sigma_0\bigcap U}\Sigma_p.
\end{equation}
Observe that $E^U_0$ is a complex analytic variety through $0$,
and that $E^{U_1}_0\subset E^{U_2}_0$ (as a germ at $0$) if
$U_2\subset U_1$. Moreover, the germ of $E^U_0$ at $0$ depends only on the
germs at $0$ of the subvarieties $\Sigma_p$ for $p\in \Sigma_0$ and, hence, does not
depend on the ball $B$. We say that $X$ is {\it essentially finite} at
$0$ if $E^U_0$ has dimension $0$ (as a germ at $0$) for every
open neighborhood $U$ of $0$.

We shall show (see Proposition
\ref{p:essvar} below) that, even if $X$ is not essentially finite at
$0$, there exists a neighborhood
$U_0$ of
$0$ such that $E^{U}_0=E^{U_0}_0$ (as germs at $0$) for every
$0\in U\subset U_0$. (This is well known in the case where $X$ is
a CR manifold.) Moreover, we shall give an alternative characterization of the stabilized germ $E^{U_0}_0$ that will be used in the proof of Theorem \ref{t:ext1}.  Let $A$ be the subvariety in $B$ defined by
\begin{equation}\Label{e:A}
A:=\{z\in B\colon \Sigma_0\subset\Sigma_z \ \text{{\rm as germs at $0$}}\}.
\end{equation}

\begin{Pro}\Label{p:essvar} Let $X\subset\bC^N$ be a
real-analytic subvariety with $0\in X$ and let $A$ be defined by \eqref{e:A}, where $B$ is a sufficiently small ball centered at $0$.  Then, there exists an open
neighborhood $U_0\subset \bC^N$ of $0$ such that $E^{U}_0=A$, as germs at $0$, for every open
$U$ with $0\in U\subset U_0$. In particular,
$E^{U}_0=E^{U_0}_0$ as germs at $0$.
\end{Pro}

\begin{proof}[Proof of Proposition $\ref{p:essvar}$]  Let $U$ be an open neighborhood of $0$ contained in $B$. Let us temporarily, in order to distinguish between germs and subvarieties, introduce the notation $D^U$ for the following representative in  $U$ of the germ at $0$ of $E^U_0$
\begin{equation}\Label{e:DU}
D^U:=\bigcap_{p\in \Sigma_0\cap U} (\Sigma_p\cap U).
\end{equation}
Observe that there exists an open neighborhood $U_0$ of $0$, contained in $B$, with the property that if $V$ is a subvariety in $B$ and $U$ is an open neighborhood of $0$ with $U\subset U_0$, then
\begin{equation}
\Sigma_0\subset V\ \text{{\rm as germs at $0$}}  \iff \Sigma_0\cap U\subset V\cap U.
\end{equation}
Indeed, it suffices to take $U_0$  so small that only the irreducible components of $\Sigma_0$ in $B$ that contain the origin meet $U_0$. As a consequence, if $U$ is an open neighborhood of $0$ with $U\subset U_0$, then the subvariety $A\cap U$ can be expressed as being those points $Z\in U$ for which $\Sigma_0\cap U\subset \Sigma_{Z}\cap U$, i.e.\
\begin{equation}
A\cap U=\{Z\in U\colon \sigma(W,\bar Z)=0\ \forall W\in U\ \text{{\rm such that }} \sigma(W,0)=0\}.
\end{equation}
On the other hand, the subvariety $D^U$, defined by \eqref{e:DU}, can be expressed in equations as follows
\begin{equation}
D^U=\{Z\in U\colon \sigma(Z,\bar W)=0\ \forall W\in U\ \text{{\rm such that }} \sigma(W,0)=0\}.
\end{equation}
Since $\sigma(Z,\bar W)=\overline{\sigma(W,\bar Z)}$, we conclude that $A\cap U=D^U$ and, hence, that $A=E^U_0$ as germs at $0$. This completes the proof of Proposition \ref{p:essvar}
\end{proof}

We shall refer to the germ at $0$ of $E^{U_0}_0$ as
the {\it essential variety} of $X$ at $0$ and denote this germ by
$E_0$. Thus, $X$ is essentially finite at $0$ if and only if $E_0=\{0\}$. When $X$
is a generic submanifold, these two notions coincide with the standard notions of
essential variety and essential finiteness (see e.g.\ \cite{BER99a}).

We shall need the
following reformulation of essential finiteness. In what follows, we shall let
$\mathcal X\subset\bC^N\times\bC^N$ denote the (local)
complexification of a real-analytic subvariety $X$ through $0$ in
$\bC^N$, i.e.\ the complex variety through $0$ in
$\bC^N\times\bC^N$ defined as the set of points $(Z,\zeta)\in \bC^N\times\bC^N$, near the origin, such that 
$$
\sigma_1(Z,\zeta)=\ldots=\sigma_d(Z,\zeta)=0
$$
where $\sigma_1(Z,\bar Z),\ldots\sigma_d(Z,\bar Z)$ are generators
of the ideal $I_\bR(X)$. For a complex analytic subvariety $V$ through $0$ in $\bC^k$, we
shall also denote by $I_{\mathcal O}(V)$ the ideal of germs at $0$ of
holomorphic functions vanishing on $V$. We shall let $I(V)$ denote the ideal
generated by $I_{\mathcal O}(V)$ in the corresponding ring of formal power series.  If $I$ is an ideal in a ring $R$, then we shall  write  $d(I)$ for the dimension of
$I$ in $R$, i.e. the dimension of the ring $R/I$ (see \cite{Ei95}; some texts,
e.g.\
\cite{AM69}, refer to the number $d(I)$ as the {\it depth} of $I$).

\begin{Lem}\Label{l:newess} Let $X$ be a real-analytic variety
through $0$ in $\bC^N$ and let $\Sigma_0$ be its Segre variety at
$0$. The following are equivalent:
\medskip

{\rm (a)} $X$ is not essentially finite at $0$.
\smallskip

{\rm (b)} There is a positive dimensional complex variety $\Gamma$
through $0$ in $\bC^N$ such that $\Gamma\times\Sigma_0^*\subset
\mathcal X$ as germs at $0$, where $\mathcal
X\subset\bC^N\times\bC^N$ denotes the complexification of $X$ and
$^*$ denotes the complex conjugate of a set (i.e.\ $S^*:=\{Z\in
\bC^N\colon \bar Z\in S\}$).
\smallskip

{\rm (c)} There is an ideal $J\subset \bC[[Z]]$ with positive dimension
such that $I(\mathcal X)\subset I(J)+I(\bC^N\times \Sigma_0^*)$.
\smallskip

{\rm (d)} There is a non-trivial formal holomorphic mapping
$\mu\colon (\bC,0)\to (\bC^N,0)$ and a neighborhood $U\subset \bC^N$
such that $\mu(\bC)\subset E^U_0$  (where $E^U_0$ is defined by
\eqref{e:essvar}), i.e.\ if $\sigma(Z,\bar Z)=(\sigma_1(Z,\bar
Z),\ldots,\sigma_d(Z,\bar Z))$ are generators of the ideal
$I_\bR(X)$, then, for every fixed $p\in \Sigma_0\cap U$,
\begin{equation}\Label{e:sigma0}
s\mapsto\sigma(\mu(s),\bar p) \ \text{{\it is identically zero as a power series.}}
\end{equation}
\end{Lem}

\begin{proof} Recall that we use the coordinates $(Z,\zeta)$ in $\bC^N\times\bC^N$.
Observe that for $p\in \bC^N$, we have
\begin{equation}\Label{e:compsig}
\begin{aligned}
\mathcal X\cap\{\zeta=\bar p\}=\Sigma_p\times\{\bar p\},\quad
\mathcal X\cap\{ Z=p\}=\{p\}\times \Sigma_p^*.
\end{aligned}
\end{equation}
Hence, $\Gamma\times\Sigma_0^*\subset \mathcal X$ means that for
every $p\in \Sigma_0$, $$\Gamma\times \{\bar p\}\subset \mathcal
X\cap\{\zeta=\bar p\}=\Sigma_p\times\{\bar p\},$$ i.e.\
$\Gamma\subset \Sigma_p$. The equivalence of (a) and (b) is a simple
consequence of this observation.

The implication (b) $\implies$ (c) is easy. Simply observe that
$\Gamma\times \Sigma_0^*=(\Gamma\times \bC^N)\cap (\bC^N\times
\Sigma_0^*)$. The conclusion of (c) now follows from (b) by taking
$J=I(\Gamma)$.

To prove the implication (c)$\implies$ (d), we let $\mu\colon
(\bC,0)\to (\bC^N,0)$ be a germ at $0$ of a non-trivial holomorphic
mapping such that $f(\mu(t))\equiv 0$ for all $f\in J$ (such exist,
by \cite{BER00}, Lemma 3.32, since $d(J)\geq 1$). It follows from
the hypothesis in (c) that there are formal power series
$d_{ij}(Z,\zeta)$ such that
\begin{equation}\Label{e:sigmahj}
\sigma_i(\mu(t),\zeta)=\sum_{j=1}^k d_{ij}(\mu(t),\zeta)h_j(\zeta),
\end{equation}
where $h_1(\zeta),\ldots,h_k(\zeta)$ generate the ideal
$I(\Sigma_0^*)\subset \bC[[\zeta]]$. Let us Taylor expand
$\sigma_i(\mu(t),\zeta)$ in $t$,
\begin{equation}
\sigma_i(\mu(t),\zeta)=\sum_{l=0}^\infty a_{il}(\zeta)t^l,
\end{equation}
and note that the coefficients $a_{il}(\zeta)$ are all holomorphic
functions of $\zeta$ in some open neighborhood $V$ of $0$. If we
Taylor expand both sides of \eqref{e:sigmahj} in $t$ and compare
coefficients, we conclude that the coefficents $a_{il}$ all belong
to the ideal $I(\Sigma_0^*)\subset \bC[[\zeta]]$. Consequently, they
vanish on $\Sigma_0^*\cap U$ for some open neighborhood $U\subset V$
of $0$. This proves (d).

To show that $(d)$ implies $(a)$, let $\mu\colon (\bC,0)\to
(\bC^N,0)$ be a non-trivial formal mapping such that
$\mu(\bC)\subset E^U_0$ for some $U$. Then $I(E^U_0)$ does not have
finite codimension (see \cite{BER00}, Lemma 3.32) and, hence, the
dimension of $E^U_0$ is positive, i.e.\ (a) holds.
\end{proof}

\section{Main lemma}

 Recall that a
formal mapping $H\colon (\bC^N,0)\to (\bC^N,0)$ is said to send (a
germ at $0$ of) a real-analytic subvariety $X\subset\bC^N$ into
another $\tilde X\subset \bC^N$, denoted $H(X)\subset \tilde X$,
if
\begin{equation}\Label{e:HXX}
\tilde \sigma\big (H(Z),\overline{H(Z)}\big)=a(Z,\bar
Z)\sigma(Z,\bar Z)
\end{equation}
holds as formal power series in $Z$ and $\bar Z$, where
$\sigma=(\sigma_1,\ldots, \sigma_d)^t$ and $\tilde \sigma=(\tilde
\sigma_1,\ldots,\tilde\sigma_{\tilde d})^t$ generate $I_\bR(X)$ and
$I_\bR(\tilde X)$, respectively, and $a(Z,\bar Z)$ is a $\tilde
d\times d$ matrix of formal power series. 

Before stating our main lemma, we review some notation and results concerning homomorphisms induced by mappings.  If $h: (\C^k_z,0) \to (\C^k_{\2z},0) $ is a formal mapping,
 then $h$ induces a homomorphism $\varphi_h: \C [[\2z ]] \to  \C[[z]] $ defined by  $\varphi_h(\tilde f)(z):=\tilde f(h(z))$.  If  $J$ is an ideal in $\C[[z]] $, then $\varphi_h^{-1}(J)$ is an ideal in $ \C [[\2z ]] $, and  $\varphi_h^{-1}(J)$ is prime if $J$ is prime.   If  $\2J$ is an ideal in $\C[[\2z]] $,  we denote by $I(\varphi_h(\2J))$ the ideal in $\C[[z]] $ generated by $\varphi_h(\2J)$. If  $\Jac h\not\equiv 0$,
where $\Jac h$ denotes the Jacobian determinant of
$h$, then \ $\varphi_h$ is injective. Indeed, in that case  any $f\in \bC\dbl \tilde Z,\tilde\zeta \dbr$
 for which $f\circ  h\equiv 0$ must be identically zero (see
e.g.\ \cite{BER99a}, Proposition 5.3.5).

We may now state
the geometric result needed to prove Theorem \ref{t:ext1}.
\begin{Lem}\Label{t:essvar} Let $X, \tilde X\subset\bC^N$ be irreducible
real-analytic subvarieties of dimension $m$ through $0$ and $\Sigma_0$, $\tilde
\Sigma_0$ their Segre varieties at $0$. Let
$H\colon(\bC^N,0)\to (\bC^N,0)$ be a finite formal mapping such that
$H(X)\subset \tilde X$. Then the following hold.
\medskip

{\rm (i)} $\varphi_H^{-1}(I(\Sigma_0))=I(\tilde\Sigma_0)$, where
$\varphi_H$ denotes the homomorphism  $\varphi_H\colon \bC[[\tilde
Z]]\to \bC [[Z]]$ induced by $H$. In particular, $\dim \tilde \Sigma_0=\dim \Sigma_0$, and if $\Sigma_0$ is irreducible, then so is $\tilde \Sigma_0$.
\smallskip

{\rm(ii)} If $\Sigma_0$  is irreducible at $0$, then $X$ is
essentially finite at $0$ if and only if $\tilde X$ is essentially
finite at $0$.
\end{Lem}

\begin{proof}[Proof of Lemma $\ref{t:essvar}$]  We begin by proving statement
(i).  We denote by $\mathcal H\colon
(\bC^N\times\bC^N,0)\to (\bC^N\times\bC^N,0)$ the complexified
formal mapping $\mathcal H(Z,\zeta)=(H(Z),\bar H(\zeta))$ and by
$\Phi=\varphi_{\mathcal H}\colon  \bC[[\tilde Z,\tilde \zeta]]\to
\bC[[Z,\zeta]]$ the induced homomorphism, i.e.\ $\Phi(\tilde
f)(Z,\zeta):=\tilde f(\mathcal H(Z,\zeta))$. The fact that $H$ sends
$X$ into $\tilde X$ can be rephrased as saying that $I(\Phi(I(\tilde
{\mathcal X})))\subset I(\mathcal X)$ or equivalently $I(\tilde
{\mathcal X})\subset \Phi^{-1}(I(\mathcal X))$, where
$\tilde{\mathcal X}$ denotes the complexification of $\tilde X$.
Since  $\Jac \mathcal H\not\equiv 0$,
\ $\Phi$ is injective.
We also claim that $\bC\dbl Z,\zeta \dbr$ is integral over
$\Phi(\bC\dbl \tilde Z,\tilde\zeta \dbr)$. To see this, note, as is
well known, that $\bC\dbl Z,\zeta\dbr$ is finitely generated over
$\Phi(\bC\dbl \tilde Z,\tilde \zeta\dbr)$ (any $f_1,\ldots,f_p\in
\bC\dbl Z,\zeta\dbr$ such that their images form a basis for the
finite dimensional vector space $\bC\dbl Z,\zeta \dbr/I(\mathcal
H_1,\ldots,\mathcal H_{2N})$ generate $\bC\dbl Z,\zeta\dbr$ over
$\Phi(\bC\dbl \tilde Z,\tilde\zeta\dbr)$). The fact that $\bC\dbl
Z,\zeta \dbr$ is integral over $\Phi(\bC\dbl \tilde
Z,\tilde\zeta\dbr)$ is now a general fact about finitely generated
modules (see e.g.\ \cite {AM69}, Proposition 5.1).  It follows that
$d(\Phi^{-1}(I(\mathcal X)))=d(I({\mathcal X}))$ (see e.g.\
\cite{Ei95},  Proposition 9.2). Moreover,
$\Phi^{-1}(I(\mathcal X))$ is a prime ideal
that contains the prime ideal $I(\tilde{\mathcal X})$. Since the
dimensions of $\mathcal X$ and $\tilde{\mathcal X}$ are the same, we
conclude that $d(\Phi^{-1}(I(\mathcal X)))=d(I(\tilde {\mathcal
X}))$ and, hence, 
\begin{equation}
\Phi^{-1}(I(\mathcal X))=I(\tilde{\mathcal X}).
\end{equation}

Observe that $\mathcal X\cap \{\zeta=0\}=\Sigma_0\times \{0\}$. By using the
specific form of the mapping $\mathcal H$, we see that
\begin{equation}
\varphi^{-1}(I(\Sigma_0))=I(\tilde \Sigma_0) \iff
\Phi^{-1}(I(\Sigma_0\times \{0\}))=I(\Sigma_0\times \{0\}),
\end{equation}
where $\varphi=\varphi_H$. By the Nullstellensatz,
\begin{equation}
I(\Sigma_0\times \{0\})=\sqrt{I(\mathcal X)+I(\zeta)}.
\end{equation}
Thus, the first identity in statement (i) is equivalent to
\begin{equation}
\sqrt{I(\tilde{\mathcal X})+I(\tilde\zeta)}=
\Phi^{-1}(\sqrt{I(\mathcal X)+I(\zeta)})
\end{equation}
We remark that the inclusion $I(\tilde{\mathcal X})\subset
\Phi^{-1}(I(\mathcal X))$ implies $I(\tilde{\mathcal
X})+I(\tilde\zeta)\subset \Phi^{-1}(I(\mathcal X)+I(\zeta))$ and
hence the inclusion
\begin{equation}
\sqrt{I(\tilde{ \mathcal X})+I(\tilde\zeta)}\subset \Phi^{-1}(\sqrt{I(\mathcal
X)+I(\zeta)})
\end{equation}
To prove the opposite inclusion, it suffices to show that $
\Phi^{-1}(I(\mathcal X)+I(\zeta))\subset\sqrt{I(\tilde {\mathcal
X})+I(\tilde\zeta)}$. Thus, we suppose that $\Phi(\tilde f)\in
I(\mathcal X)+I(\zeta)$ and must prove $\tilde f^k\in I(\tilde
{\mathcal X})+I(\tilde\zeta)$ for some $k$. Since $\bar H$ is a
formal finite mapping, there is an integer $m$ such that
$I(\zeta)^m\subset I(\Phi(I(\tilde\zeta)))$. We conclude that, for
any $k\geq m$,
\begin{equation}
\Phi(\tilde f)^k=\Phi(\tilde f^k)\in I(\mathcal
X)+I(\zeta)^k\subset I(\mathcal X)+I(\Phi(I(\tilde\zeta))).
\end{equation}
Hence, we have
\begin{equation}\Label{e:expan1}
\Phi(\tilde f^k)=g_k(Z,\zeta)+\sum_{j=1}^N a_{jk}(Z,\zeta)\bar H_j(\zeta),
\end{equation}
where $g_k\in I(\mathcal X)$ and $a_{jk}\in \bC[[Z,\zeta]]$. Let
$(\alpha_j,\beta_j)$, for $j=0,\ldots,p$, be multi-indices such that
the images of $Z^{\alpha_j}\zeta^{\beta_j}$ generate the finite
dimensional vector space $\bC[[Z,\zeta]]/I(\mathcal H)$. Observe
that $(0,0)$ is necessarily one of these; we order the multi-indices
in such a way that $(\alpha_0,\beta_0)=(0,0)$. We may then, as is easy to verify, write
each $a_j(Z,\zeta)$ in the following way  
\begin{equation}\Label{e:new}
a_{jk}(Z,\zeta)=\sum_{l=0}^p b_{jlk}(\mathcal
H(Z,\zeta))Z^{\alpha_l}\zeta^{\beta_l},
\end{equation}
where $b_{00k}(0)=0$ (since $\Phi(\tilde f^k)(0)=0$). By substituting \eqref{e:new} in \eqref{e:expan1}, it follows
 that
\begin{equation}\Label{e:expan2}
\Phi(\tilde f^k)=g_k(Z,\zeta)+\sum_{l=0}^p\sum_{j=1}^N
b_{jlk}(\mathcal H(Z,\zeta))\bar
H_j(\zeta)Z^{\alpha_l}\zeta^{\beta_l}.
\end{equation}
Let $S$ and $\tilde S$ denote rings $\bC[[Z,\zeta]]/I(\mathcal X)$
and $\bC[[\tilde Z,\tilde \zeta]]/I(\tilde{\mathcal X})$,
respectively. We shall denote by $f^*$ and $A^*$ the images of any
element $f\in \bC[[Z,\zeta]]$ and subset $A\subset\bC[[Z,\zeta]]$ in
$S$, and similarly for images in $\tilde S$. Since
$\Phi(I(\tilde{\mathcal X})\subset I({\mathcal X})$, the
homomorphism $\Phi$ induces a homomorphism $\Phi^*\colon \tilde S\to
S$. The facts that $\Phi$ is injective and $\Phi^{-1}(I( \mathcal
X))=I(\tilde{\mathcal X})$ imply that $\Phi^*$ is injective. Let
$\tilde R$ denote the ring with identity generated by
$I(\tilde\zeta)^*\subset\tilde S$, i.e.\ $$\tilde R:=\{\tilde
g^*+c^*\colon \tilde g\in I(\tilde\zeta),\ c\in \bC\}.$$ Denote by
$R$ the ring $\Phi^*(\tilde R)\subset S$ and let $\mathcal N$ denote the
$R$-module generated by $(Z^\alpha_j\zeta^{\beta_j})^*$ for
$j=0,\ldots,p$. Clearly, $\mathcal N$ is not annihilated by any elements of
$S$. Let $s$ denote $\Phi(\tilde f^m)^*\in S$ and observe that
\begin{equation}\Label{e:expan3} s=\Phi(\tilde
f^m)^*=\sum_{l=0}^p\sum_{j=1}^N \big(b_{jlm}(\mathcal
H(Z,\zeta))\bar H_j(\zeta)Z^{\alpha_l}\zeta^{\beta_l}\big)^*.
\end{equation}
We claim that $s\mathcal N\subset \mathcal N$. Indeed, if $f^*\in \mathcal N$, then
$$
f^*=\sum_{i=0}^p \big(d_{i}(\mathcal
H(Z,\zeta))Z^{\alpha_i}\zeta^{\beta_i}\big)^*,
$$
where $d_{l}=\tilde g+c$ for some $\tilde g\in I(\tilde\zeta)$ and
$c\in \bC$, and hence
\begin{equation}\Label{e:mult} sf^*=
\sum_{l=0}^p\sum_{i=0}^p\sum_{j=1}^N \big(b_{jlm}(\mathcal
H(Z,\zeta))d_{i}(\mathcal H(Z,\zeta))\bar
H_j(\zeta)Z^{\alpha_l+\alpha_i}\zeta^{\beta_l+\beta_i}\big)^*.
\end{equation}
Any monomial $Z^\alpha\zeta^\beta$ can be written
\begin{equation}\Label{e:mono}
Z^\alpha\zeta^\beta=\sum_{q=0}^p e_{\alpha\beta q}(\mathcal
H(Z,\zeta))Z^{\alpha_q}\zeta^{\beta_q}.
\end{equation}
By substituting \eqref{e:mono} in \eqref{e:mult}, we conclude that
$sf^*\in\mathcal N$, as claimed. It follows from \cite{Ei95}, Corollary 4.6,
that $s$ is integral over $R$, i.e.\
\begin{equation}\Label{e:integral}
s^r+\sum_{k=0}^{r-1} (h_k(\mathcal H(Z,\zeta)))^* s^k=0,
\end{equation}
where $h_k=\tilde g_k+c_k$ for $\tilde g_k\in I(\tilde \zeta)$ and
$c_k\in \bC$. Since $\Phi^*$ is injective (and $s=\Phi(\tilde
f^m)^*$), we conclude that
\begin{equation}\Label{e:integral2}
(\tilde f^{rm}(\tilde Z,\tilde\zeta))^*+\sum_{k=0}^{r-1} (h_k(\tilde
Z,\tilde \zeta))^* (\tilde f^{km}(\tilde Z,\tilde\zeta))^*=0.
\end{equation}
Since $h_k=\tilde g_k+c_k$, we can rewrite this as
\begin{equation}\Label{e:integral3}
(\tilde f^{rm}(\tilde Z,\tilde\zeta))^*+\sum_{k=0}^{r-1} (g_k(\tilde
Z,\tilde \zeta))^* (\tilde f^{km}(\tilde Z,\tilde\zeta))^* +
\sum_{k=1}^{r-1} c_k^* (\tilde f^{km}(\tilde
Z,\tilde\zeta))^*=-c_0^*.
\end{equation}
First, observe that the left hand side of \eqref{e:integral3} is in
the maximal ideal of $\tilde S$ and, hence, $c_0^*=0$. Let us write
$c_r^*=1$, and let $k_0$ be the smallest integer in $\{1,\ldots,r\}$
such that $c_{k_0}^*\neq 0$. We may then rewrite \eqref{e:integral3}
as
\begin{equation}\Label{e:integral4}
(\tilde f^{mk_0}(\tilde Z,\tilde\zeta))^*\bigg(\sum_{k=k_0}^{r}
c_k^* (\tilde f^{m(k-k_0)}(\tilde Z,\tilde\zeta))^*\bigg)
=-\sum_{k=0}^{r-1} (g_k(\tilde Z,\tilde \zeta))^* (\tilde
f^{km}(\tilde Z,\tilde\zeta))^*.
\end{equation}
The right hand side of \eqref{e:integral4} is in $I(\tilde\zeta)^*$
and, hence, so is then the left hand side. Observe that
$$
\tilde f^{m(r-k_0)}(\tilde Z,\tilde\zeta))^* + \sum_{k=k_0}^{r-1}
c_k^* (\tilde f^{m(k-k_0)}(\tilde Z,\tilde\zeta))^*
$$ is a unit, since $c_{k_0}^*\neq 0$. We conclude that $(\tilde
f^{mk_0})^*\in I(\tilde\zeta)^*$ and, hence, $\tilde f^{mk_0}\in
I(\tilde{\mathcal X})+I(\tilde\zeta)$, as desired. This completes
the proof of the identity  $\varphi_H^{-1}(I(\Sigma_0))=I(\tilde\Sigma_0)$.  We observe that if $I(\Sigma_0)$ is prime, then so is $I(\tilde\Sigma_0)$.  The fact that  $d(I(\tilde \Sigma_0))=d(I(\Sigma_0))$ (or equivalently $\dim \Sigma_0=\dim\tilde  \Sigma_0$)  follows from the assumption that $H$ is a finite mapping (cf.\ the use of
\cite{Ei95}, Proposition 9.2 above).  This completes
the proof of (i) in Lemma
\ref{t:essvar}.

\section{Proof of {\rm (ii)} of Lemma \ref{t:essvar}}

We now proceed with the proof of statement (ii) in Lemma
\ref{t:essvar}. The fact that $\Sigma_0$ is assumed to be irreducible at $0$ implies that $I(\Sigma_0)$ is a prime ideal. By part
(i) of Lemma \ref{t:essvar}, it follows that
$I(\tilde\Sigma_0)(=\varphi_H^{-1}(I(\Sigma_0)))$ is a prime ideal and 
$d(\varphi_H^{-1}(I(\Sigma_0)))=d(I(\Sigma_0))$.  Equivalently,
$\tilde\Sigma_0$ is an irreducible complex analytic variety of the
same dimension as $\Sigma_0$.

We begin by proving that if $X$ is essentially finite at $0$, then
$\tilde X$ is also essentially finite at $0$. We shall show the
logical negation of this statement. Thus, we shall assume that
$\tilde X$ is not essentially finite at $0$. By Lemma
\ref{l:newess}, there exists a positive dimensional complex analytic
variety $\tilde \Gamma\subset \bC^N$ through $0$ such that $\tilde
\Gamma\times \tilde \Sigma_0^*\subset \tilde{\mathcal X}$. We may
assume that $\tilde \Gamma$ is irreducible.

  \begin{Lem} \Label{l:step1} Let $\tilde \Gamma$ be as above and let $\tilde
{\frak P}$ be the (prime) ideal in $\bC\dbl \tilde Z,\tilde
\zeta\dbr$ generated by the ideal of $\Gamma\times \tilde\Sigma_0^*$
(in the ring $\bC\{\tilde Z,\tilde\zeta \}$ of convergent power
series).  Let $\Phi:=\varphi_{\mathcal H}$ be the homomorphism defined in the beginning of the proof of Lemma $\ref{t:essvar}$. 
If
\begin{equation}\Label{e:primary}
I(\Phi(\tilde{\frak P}))=\bigcap_{j=1}^m P_j
\end{equation}
is a primary decomposition of $I(\Phi(\tilde{\frak P}))$, then there exists $j_0$ such that
\begin{equation} I(\mathcal X) \subset \sqrt{P_{j_0}\ } \ {\rm{and}} \ \ I(P_{j_0})=d(\tilde {\frak P}).
\end{equation}
\end{Lem}

\begin{proof}  Suppose, in order to reach a contradiction, that this
is not the case. Then, each of the ideals $I(\mathcal X)+P_j$ must
satisfy $$d(I(\mathcal X)+P_j)\leq d(P_j)-1\leq d(\tilde{\frak
P})-1.$$ It follows that
\begin{equation}
d\bigg(\bigcap_{i=1}^m(I(\mathcal X)+P_j)\bigg)\leq d(\tilde{\frak
p})-1.
\end{equation}
Consider the induced homomorphism
\begin{equation}\Label{e:Phi*}
\Phi^*\colon \bC\dbl \tilde Z,\tilde \zeta\dbr/I(\tilde {\mathcal
X})\to \bC\dbl Z, \zeta\dbr/I(\mathcal X).
\end{equation}
If we, as above, use $J^*$ to denote the image of an ideal $J\subset
\bC\dbl\tilde Z,\tilde\zeta\dbr$ in $\bC\dbl \tilde Z,\tilde
\zeta\dbr/I(\tilde {\mathcal X})$ and similarly for images in
$\bC\dbl Z, \zeta\dbr/I(\mathcal X)$, then we have
\begin{equation}
I(\Phi^*({\tilde {\frak P}}^*))=\bigg(\bigcap_{j=1}^m P_j\bigg)^*.
\end{equation}
Observe that
$$
\bigg(\bigcap_{j=1}^m P_j^*\bigg)^m\subset \bigg(\bigcap_{j=1}^m
P_j\bigg)^* \subset \bigcap_{j=1}^m P_j^*.$$ Since
$P_j^*=(P_j+I(\mathcal X))^*$ and, hence, $d(P^*_j)\leq
d(\tilde{\frak  P}^*)-1$, we conclude that
\begin{equation}\Label{e:nottrue}
d(I(\Phi^*({\tilde {\frak P}}^*)))\leq d(\tilde{\frak  P}^*)-1.
\end{equation}
We shall show that this is a contradiction. Recall that the
homomorphism $\Phi^*$ is injective (see the proof of (i) above).
Clearly, $\bC\dbl Z, \zeta\dbr/I(\mathcal X)$ is integral over
$\Phi^*(\bC\dbl \tilde Z,\tilde \zeta\dbr/I(\tilde {\mathcal X}))$,
since $\bC\dbl Z,\zeta\dbr$ is integral over $\Phi(\bC\dbl\tilde
Z,\tilde\zeta\dbr)$. Now, the fact that \eqref{e:nottrue} cannot
hold is an immediate consequence of  Lemma \ref{l:lemma1} below.\end{proof}

The proof of Lemma \ref{l:step1} follows, as mentioned above, from the following commutative algebra lemma.
For the reader's convenience, we include a
proof; we have been unable to find an exact reference for this
statement.

\begin{Lem}\Label{l:lemma1} Let $A$ and $B$ be (commutative)
rings, and $\psi\colon A\to B$ an injective homomorphism such that
$B$ is integral over $\psi(A)$. Then, for any ideal $J\subset A$,
$$
d(I(\psi(J)))= d(J).
$$
\end{Lem}

\begin{proof}[Proof of Lemma $\ref{l:lemma1}$] Since $\psi$ is injective,
we may identify $A$ with the subring $\psi(A)\subset B$. Thus, $J$
is identified with its image $\psi(J)$ and $I(J)$ is the ideal in
$B$ generated by $J$. Let
\begin{equation}\Label{e:dpsiJ}
I(J)\subset\frak q_0\subsetneq\frak q_1\subsetneq\ldots
\subsetneq\frak q_r,
\end{equation}
be a chain of prime ideals. Then $\frak p_i=\frak q_i\bigcap A$
($=\psi^{-1}(\frak q_i)$) are prime ideals and since $J\subset
I(J)\bigcap A$, we obtain
\begin{equation}\Label{e:hpsiJ}
J\subset\frak p_0\subset\frak p_1\subset\ldots \subset\frak p_r.
\end{equation}
By Corollary 4.18 in \cite{Ei95}, we also have strict inclusions
$\frak p_i\subsetneq \frak p_{i+1}$ for $i=0,\ldots, r-1$. It
follows that $d(J)\geq d(I(J))$.

To prove the opposite inequality, let
\begin{equation}\Label{e:dJ}
J\subset\frak p'_0\subsetneq\frak p'_1\subsetneq\ldots
\subsetneq\frak p'_s,
\end{equation}
be a chain of prime ideals. By the Going Up Lemma (see
\cite{Ei95}, Proposition 4.15), there is a prime ideal $\frak
q'_0$ in $B$ such that $\frak p_0'=\frak q'_0\bigcap A$. Hence,
$I(J)\subset \frak q'_0$. By inductively applying the Going Up
Lemma, we find prime ideals $\frak q'_1,\ldots,\frak q'_s$ such
that $\frak p'_i=\frak q_i\bigcap A$ and\
\begin{equation}\Label{e:dpsiJ}
I(J)\subset\frak q'_0\subset\frak q'_1\subset\ldots \subset\frak
q'_s.
\end{equation}
Clearly, we have strict inclusions $\frak q'_i\subsetneq \frak
q'_{i+1}$ for $i=0,\ldots, s-1$ (since $\frak p'_i=\frak
q_i\bigcap A$ and the $\frak p_i$ are distinct). This proves the
opposite inequality $d(I(J))\geq d(j)$, which completes the proof
of the lemma.
\end{proof}

To complete the proof of statement (ii) in Lemma \ref{t:essvar}, we shall need the following lemma.

\begin{Lem}\Label{l:step2} Let $\varphi:=\varphi_H\colon \bC \dbl
\tilde Z\dbr\to \bC\dbl Z\dbr$ denote the homomorphism induced by
the formal mapping $H\colon (\bC^N,0)\to (\bC^N,0)$ and
$\psi:=\varphi_{\bar H}\colon \bC \dbl \tilde \zeta\dbr\to \bC\dbl \zeta\dbr$ the homomorphism induced by $\bar H\colon (\bC^N,0)\to
(\bC^N,0)$. Every minimal prime $\sqrt{P}=\sqrt{P_j}$ in the factorization \eqref{e:primary}  is
of the form
\begin{equation}\Label{e:claim2}
\sqrt{P}=I(\frak p)+I(\frak q),
\end{equation}
where $\frak p$ is a minimal prime of $I(\varphi(I(\tilde
\Gamma)))\subset \bC\dbl Z\dbr$ and $\frak q$ is a minimal prime
of $I(\psi(I(\tilde\Sigma^*_0)))\subset \bC\dbl\zeta\dbr$.
\end{Lem}

\begin{proof}  We first observe that if $\frak p\subset \bC\dbl Z\dbr$ and
$\frak q\subset \bC\dbl \zeta\dbr$ are prime ideals, then $I(\frak
p)+I(\frak q)\subset \bC\dbl Z,\zeta \dbr$ is also prime. For,
$\bC\dbl Z,\zeta\dbr=\bC\dbl Z\dbr\otimes \bC\dbl \zeta\dbr$
(where the tensor product is over $\bC$) and, by Theorem III.14.35
of \cite{ZS58},
\begin{equation}
\bC\dbl Z,\zeta\dbr/(I(\frak p)+I(\frak q))= (\bC\dbl Z\dbr/\frak
p)\otimes (\bC\dbl \zeta\dbr/\frak q).
\end{equation}
Since $\frak p$ and $\frak q$ are prime, $\bC\dbl Z\dbr/\frak p$ and
$\bC\dbl \zeta\dbr/\frak q$ are integral domains. It follows that
$\bC\dbl Z\dbr/\frak p)\otimes (\bC\dbl \zeta\dbr/\frak q)$ is an
integral domain, since $(\bC\dbl Z\dbr/\frak p)\otimes (\bC\dbl
\zeta\dbr/\frak q)\subset K\otimes K'$, where $K$ and $K'$ denote
the quotient fields of $\bC\dbl Z\dbr/\frak p$ and $\bC\dbl
\zeta\dbr/\frak q$ respectively, and $K\otimes K'$ is an integral
domain by Corollary III.15.1 of \cite{ZS58}. This proves that
$I(\frak p)+I(\frak q)$ is prime. By considering maximal chains of
prime ideals containing $I(\frak p)+I(\frak q)$ constructed in an
obvious way from maximal chains of prime ideals containing $\frak p$
and $\frak q$, we also deduce that $d(I(\frak p)+I(\frak q))=d(\frak
p)+d(\frak q)$. Now, let $\frak p_1,\ldots, \frak p_k$ and $\frak
q_1,\ldots, \frak q_l$ be the minimal primes of $I(\varphi(I(\tilde
\Gamma)))$ and $I(\psi(I(\tilde\Sigma^*_0)))$, respectively, i.e.\
$d(\frak p_i)=d(I(\varphi(I(\tilde \Gamma))))$, $d(\frak
q_j)=d(I(\psi(I(\tilde\Sigma^*_0))))$ and
$$
\sqrt{I(\varphi(I(\tilde \Gamma)))}=\bigcap _{i=1}^k \frak
p_i,\quad \sqrt{I(\psi(I(\tilde\Sigma^*_0)))}=\bigcap _{i=j}^l
\frak q_j.
$$
To prove the decomposition \eqref{e:claim2}, we shall show that
\begin{equation}\Label{e:P=p+q}
\sqrt{I(\Phi(\tilde{\frak P}))}=\bigcap _{i=1}^k \bigcap _{j=1}^l
(I(\frak p_i)+I(\frak q_j)).
\end{equation}
The uniqueness of the minimal primes of an ideal then implies that
\eqref{e:claim2} must hold. The inclusion
$$
\sqrt{I(\Phi(\tilde{\frak P}))}\subset\bigcap _{i=1}^k \bigcap
_{j=1}^l (I(\frak p_i)+I(\frak q_j))
$$ is easy to prove and the details of this are left to the
reader. Now, let $h\in \cap _{i,j}(I(\frak p_i)+I(\frak q_j))$.
Fix a $j\in \{1,\ldots, l\}$ and let $f_i\in I(\frak p_i)$,
$g_i\in I(\frak q_j)$ such that $h=f_i+g_i$ for $i=1,\ldots, k$.
It follows that
$$
h^k=\prod_{i=1}^k (f_i+g_i).
$$
Since $f_{i_1}\ldots f_{i_{k-r}}g_{i_{k-r+1}}\ldots g_{i_k}$
belongs to $I(\frak q_j)$ whenever $r\geq 1$ and $f_1\ldots f_k$
belongs to $\cap_i I(\frak p_i)$, we conclude that
\begin{equation}
h^k= f'_j+g'_j,
\end{equation}
where $ f'_j\in\cap_i I(\frak p_i)$ and $g_j'\in I(\frak q_j)$ for
every $j=1,\ldots,l$. A similar argument shows that $h^{kl}=f+g$
with $f\in \cap_i I(\frak p_i)$ and $g\in \cap_j I(\frak q_j)$. We
conclude that
$$
h^{kl}\in \cap_i I(\frak p_i)+\cap_j I(\frak q_j)
$$
or, equivalently,
\begin{equation}
h\in \sqrt{\sqrt{I(\varphi(I(\tilde
\Gamma)))}+\sqrt{I(\psi(I(\tilde\Sigma^*_0)))}}
=\sqrt{I(\Phi(\tilde{\frak P}))}.
\end{equation}
This proves \eqref{e:P=p+q} and, hence, also the lemma.
\end{proof}

We now return to the proof of statement (ii) in Lemma
\ref{t:essvar}. Let $P:=P_{j_0}$ be a primary ideal in
\eqref{e:primary} as given by Lemma \ref{l:step1}. Let $\frak p$ and
$\frak q$ be as in Lemma \ref{l:step2} such that \eqref{e:claim2} is
satisfied. In particular, $d(\frak p)\geq 1$. By evaluating at
$Z=0$, we deduce that $I(\Sigma^*_0)\subset \frak q$. Since
$$d(\frak q)=d(I(\tilde \Sigma^*_0))= d(I(\Sigma^*_0)),$$ where the
latter identity follows from statement (i) in Lemma \ref{t:essvar},
and both $\frak q$ and $I(\Sigma^*_0)$ are primes, we conclude that,
in fact, $\frak q=I(\Sigma^*_0)$. Hence, $I(\mathcal X)\subset
I(\frak p)+I(\bC^N\times \Sigma_0^*)$. The fact that $\mathcal X$ is
not essentially finite at $0$ now follows from Lemma \ref{l:newess},
part (c) with $J=\frak p$. This completes the proof of the
implication ``$X$ is essentially finite at $0$" $\implies$ ``$\tilde
X$ is essentially finite at $0$".

To finish the proof of (ii), we must show the converse implication,
namely ``$\tilde X$ is essentially finite at $0$" $\implies$ ``$X$
is essentially finite at $0$". Again, we shall prove the logical
negation of this statement.  Thus, suppose that $X$ is not
essentially finite at $0$ and let $\Gamma$ be an irreducible complex
analytic variety through $0$ in $\bC^N$ such that
$\Gamma\times\Sigma^*_0\subset I(\mathcal X)$. Let $\frak
p:=I(\Gamma)$, $\frak q:=I(\Sigma^*_0)$, and observe, as above, that
$\frak P:=I(\frak p)+I(\frak q)$ is a prime ideal such that
$I(\mathcal X)\subset \frak P$. Let $\tilde {\frak P}$ be the prime
ideal $\Phi^{-1}(\frak P)$ and observe that $I(\tilde{\mathcal
X})\subset \tilde{\frak P}$. Recall the homomorphisms $\varphi\colon
\bC[[\tilde Z]]\to \bC[[Z]]$ and $\psi\colon \bC[[\tilde \zeta]]\to
\bC[[\zeta]]$ induced by $H(z)$ and $\bar H(\zeta)$. We claim that
\begin{equation}\Label{e:pid} \tilde{\frak P}=I(\phi^{-1}(\frak
p))+I(\psi^{-1}(\frak q)).
\end{equation}
 Indeed, we have
$\Phi^{-1}(I(\frak p))+\Phi^{-1}(I(\frak q))\subset \Phi^{-1}(\frak
P)$, $\Phi^{-1}(I(\frak p))=I(\phi^{-1}(\frak p))$, and
$\Phi^{-1}(I(\frak q))=I(\psi^{-1}(\frak q))$, which together imply
the inclusion
\begin{equation}\Label{e:pinc}
I(\phi^{-1}(\frak p))+I(\psi^{-1}(\frak q))\subset \tilde{\frak P}.
\end{equation}
Both sides of \eqref{e:pinc} are prime ideals and the following
chain of identities of dimensions follow from the results above
$$d(\tilde{\frak P})=d(\frak P)=d(\frak p)+d(\frak
q)=d(\phi^{-1}(\frak p))+d(\psi^{-1}(\frak q))=d(I(\phi^{-1}(\frak
p))+I(\psi^{-1}(\frak q))).$$ This implies the desired identity
\eqref{e:pid}. By (i) of the lemma, we have $\psi^{-1}(\frak
q))=I(\tilde\Sigma^*_0)$. Since $d(\phi^{-1}(\frak p))\geq 1$, an
argument analogous to that used to complete the proof of the
implication ``$X$ is essentially finite at $0$" $\implies$ ``$\tilde
X$ is essentially finite at $0$" above now shows that $\tilde X$ is
not essentially finite at $0$. This completes the proof of Lemma
\ref{t:essvar}.
\end{proof}

\section{Proofs of Theorems \ref{t:ext1} and \ref{t:ext2}}
\begin{proof}[Proof of Theorem  $\ref{t:ext1}$] We first prove the implication
(i)$\implies$(ii) in Theorem \ref{t:ext1}.  Let $H\colon (\bC^N,0)\to (\bC^N,0)$ be the formal holomorphic mapping associated to
$f$ sending $M$ into $\tilde X$, and $\phi_H\colon \bC[[\tilde Z]]\to \bC[[Z]]$ the homomorphism induced by $H$. We shall use the following sufficient criterion from \cite{MMZ02} (see Theorem 1.1 and Remark 1.2 in \cite{MMZ02}) for $f$ to be real-analytic at $0$ (or, equivalently, extend holomorphically near $0$): $f$ is real-analytic near $0$ if the germ at $0$ of the variety
\begin{equation}\Label{e:formalMMZ}
C:=\{\tilde Z\colon \phi_H(I(\tilde \Sigma_{\tilde Z}))\subset  I(\Sigma_0)\}
\end{equation}
reduces to the single point $\{0\}$.  We should perhaps point out that if $0$ does not belong to $\tilde \Sigma_{\tilde Z}$, then $I(\tilde\Sigma_{\tilde Z})$ is the whole ring $\bC[[\tilde Z]]$ (and, hence, $\tilde Z$ does not belong to $C$). Observe that $C$ can also be expressed as
\begin{equation}\Label{e:formalMMZ1}
C:=\{\tilde Z\colon I(\tilde \Sigma_{\tilde Z}))\subset \phi_H^{-1}( I(\Sigma_0))\}
\end{equation}
By using Lemma \ref{t:essvar} (i), we observe that $ I(\tilde \Sigma_{\tilde Z}))\subset \phi_H^{-1}( I(\tilde \Sigma_0))$ implies that $ I(\tilde \Sigma_{\tilde Z}))\subset I(\tilde\Sigma_0)$, i.e.\  $\tilde \Sigma_0\subset \tilde \Sigma_{\tilde Z}$. It follows that $C\subset \tilde A$ as germs at $0$, where $\tilde A$ is given by \eqref{e:A} using $\tilde X$ instead of $X$ .  Since $M$ is assumed to be essentially finite at $0$,  it follows that $\tilde X$ is also essentially finite at $0$ by Lemma \ref{t:essvar} (ii). By Proposition \ref{p:essvar}, the germ of $\tilde A$ at $0$, and hence that of $C$, reduces to the point $\{0\}$. The real-analyticity of $f$ near $0$ now follows from the above mentioned result from \cite{MMZ02}.

To prove (ii)$\implies$(i), we let $H\colon (\bC^N,0)\to (\bC^N,0)$ be the finite holomorphic mapping extending $f$ near $0$. We further let  $\mathcal M\subset \bC^N\times \bC^N$ denote the local complexification of $M$ near $0$ and $\mathcal H\colon (\bC^N\times\bC^N,0)\to (\bC^N\times\bC^N,0)$ the complexification of the mapping $H$. Since $\mathcal H$ is a finite holomorphic mapping and $\mathcal M$ a complex submanifold, we conclude, by Remmert's proper mapping theorem, that the image $\tilde{\mathcal X}:=\mathcal H(\mathcal M)$ is an irreducible complex analytic subvariety of the same (complex) dimension as $\mathcal M$.  Let $\tilde X$ denote the real-analytic subvariety of $\bC^N$ obtained by intersecting $\tilde{\mathcal X}$ with the anti-diagonal $\{(\tilde Z,\tilde \zeta)\colon \tilde\zeta=\bar{\tilde Z}\}$. It is easy to check that $H$ maps $M$ into $\tilde X$. Moreover, since $\det (\partial H/\partial Z)$ does not vanish identically on $M$, the dimension of $\tilde X$ is at least that of $M$. On the other hand,
since $\tilde X$ is the intersection between the totally real manifold $\{(\tilde Z,\tilde \zeta)\colon \tilde\zeta=\bar{\tilde Z}\}$ and $\tilde{\mathcal X}$, its dimension is at most equal to the complex dimension of $\tilde{\mathcal X}$. Since
\begin{equation}
\dim_\bC \tilde{\mathcal X}=\dim_\bC \mathcal M=\dim_\bR M,
\end{equation}
we conclude that $\dim_\bR \tilde X=\dim_\bR M$. The fact that $\tilde X$ is irreducible now follows from the fact that $\tilde{\mathcal X}$ is its complexification and $\tilde{\mathcal X}$ is irreducible. This completes the proof of Theorem \ref{t:ext1}.
\end{proof}

\begin{proof}[Proof of Theorem  $\ref{t:ext2}$]  The proof of (ii) $\implies$ (i) in Theorem  $\ref{t:ext2}$ is exactly the same as the proof of (ii) $\implies$ (i) in Theorem  $\ref{t:ext1}$.  The proof of (i) $\implies$ (ii) in Theorem  $\ref{t:ext2}$ is follows the same reasoning as that of the proof of (i) $\implies$ (ii) in Theorem  $\ref{t:ext1}$, except that we apply a result from \cite{Mirg} instead of the result from \cite{MMZ02} used in the proof above. Indeed, Theorem 9.1 in \cite{Mirg} combined with the remark preceding it shows that the formal mapping $H$ in Theorem \ref{t:ext2} is convergent if the variety $C$ given by \eqref{e:formalMMZ} reduces to the single point $\{0\}$. This completes the proof of Theorem \ref{t:ext2}.
\end{proof}

\end{document}